\documentclass[a4paper,12pt]{amsart}
\usepackage{amsfonts}
\usepackage{amssymb}
\usepackage{ifthen}
\usepackage{amscd}
\usepackage{amsxtra}
\usepackage{graphicx}
\usepackage{color}
\nonstopmode \numberwithin{equation}{section}
\newtheorem{thm}{Theorem}
\newtheorem{lem}{Lemma}
\newtheorem{cor}{Corollary}
\newtheorem{prop}{Proposition}

\newtheorem{cl}{Claim}
\newtheorem{ca}{Case}
\newtheorem{sca}{Subcase}
\newtheorem{scl}{Subclaim}
\newtheorem{conj}{Conjecture}

\theoremstyle{definition}
\newtheorem{defn}{Definition}

\newtheorem{op}[equation]{Open Problem}
\newtheorem{ques}[equation]{Question}
\newtheorem{rem}{Remark}[section]
\newtheorem{exam}[equation]{Example}

\newcounter {own}
\def\theown {\thesection       .\arabic{own}}

\newenvironment{pf}[1][]{%
 \vskip 3mm
 \noindent
 \ifthenelse{\equal{#1}{}}%
  {{\slshape Proof. }}%
  {{\slshape #1.} }%
 }%
{\qed\bigskip}

\newcounter{alphabet}
\newcounter{tmp}
\newenvironment{Thm}[1][]{\refstepcounter{alphabet}%
\bigskip%
\noindent%
{\bf Theorem \Alph{alphabet}}%
\ifthenelse{\equal{#1}{}}{}{ (#1)}%
{\bf .} \itshape}{\vskip 8pt}

\makeatletter
\newcommand{\Ref}[1]{\@ifundefined{r@#1}{}{\setcounter{tmp}{\ref{#1}}\Alph{tmp}}}
\makeatother

\newenvironment{Lem}[1][]{\refstepcounter{alphabet}%
\bigskip%
\noindent%
{\bf Lemma \Alph{alphabet}}%
{\bf .} \itshape}{\vskip 8pt}

\newcommand{\IC}{{\mathbb C}}
\newcommand{\ID}{{\mathbb D}}

\newcommand{\diam}{{\operatorname{diam}}}

\def\be{\begin{equation}}
\def\ee{\end{equation}}

\newcommand{\bee}{\begin{enumerate}}
\newcommand{\eee}{\end{enumerate}}

\newcommand{\blem}{\begin{lem}}
\newcommand{\elem}{\end{lem}}
\newcommand{\bthm}{\begin{thm}}
\newcommand{\ethm}{\end{thm}}
\newcommand{\bcor}{\begin{cor}}
\newcommand{\ecor}{\end{cor}}
\newcommand{\beg}{\begin{exam}}
\newcommand{\eeg}{\end{exam}}
\newcommand{\begs}{\begin{examples}}
\newcommand{\eegs}{\end{examples}}
\newcommand{\bdefe}{\begin{defn}}
\newcommand{\edefe}{\end{defn}}
\newcommand{\bprob}{\begin{prob}}
\newcommand{\eprob}{\end{prob}}
\newcommand{\bques}{\begin{ques}}
\newcommand{\eques}{\end{ques}}
\newcommand{\bei}{\begin{itemize}}
\newcommand{\eei}{\end{itemize}}
\newcommand{\bcon}{\begin{conj}}
\newcommand{\econ}{\end{conj}}
\newcommand{\bop}{\begin{op}}
\newcommand{\eop}{\end{op}}

\newcommand{\bca}{\begin{ca}}
\newcommand{\eca}{\end{ca}}
\newcommand{\bsca}{\begin{sca}}
\newcommand{\esca}{\end{sca}}

\newcommand{\bcl}{\begin{cl}}
\newcommand{\ecl}{\end{cl}}

\newcommand{\bscl}{\begin{scl}}
\newcommand{\escl}{\end{scl}}

\newcommand{\bcons}{\begin{conjs}}
\newcommand{\econs}{\end{conjs}}
\newcommand{\bprop}{\begin{propo}}
\newcommand{\eprop}{\end{propo}}
\newcommand{\br}{\begin{rem}}
\newcommand{\er}{\end{rem}}
\newcommand{\brs}{\begin{rems}}
\newcommand{\ers}{\end{rems}}
\newcommand{\bo}{\begin{obser}}
\newcommand{\eo}{\end{obser}}
\newcommand{\bos}{\begin{obsers}}
\newcommand{\eos}{\end{obsers}}
\newcommand{\bpf}{\begin{pf}}
\newcommand{\epf}{\end{pf}}
\newcommand{\ba}{\begin{array}}
\newcommand{\ea}{\end{array}}
\newcommand{\beq}{\begin{eqnarray}}
\newcommand{\beqq}{\begin{eqnarray*}}
\newcommand{\eeq}{\end{eqnarray}}
\newcommand{\eeqq}{\end{eqnarray*}}

\newcommand{\ds}{\displaystyle}

\newcounter{minutes}
\divide\time by 60
\newcounter{hours}
\multiply\time by 60 \addtocounter{minutes}{-\time}

\begin{document}

\bibliographystyle{amsplain}
\title []
{linear measure and  $K$-quasiconformal  harmonic mappings}

\def\thefootnote{}
\footnotetext{ \texttt{\tiny File:~\jobname .tex,
          printed: \number\day-\number\month-\number\year,
          \thehours.\ifnum\theminutes<10{0}\fi\theminutes}
} \makeatletter\def\thefootnote{\@arabic\c@footnote}\makeatother

\author{Shaolin Chen}
 \address{ College of Mathematics and
Statistics, Hengyang Normal University, Hengyang, Hunan 421008,
People's Republic of China.} \email{mathechen@126.com}

\author{Gang  Liu}
\address{G. Liu, College of Mathematics and
Statistics, Hengyang Normal University, Hengyang, Hunan 421008,
People's Republic of China.} \email{liugangmath@sina.cn}

\author{Saminathan Ponnusamy}
\address{S. Ponnusamy,
Indian Statistical Institute (ISI), Chennai Centre, SETS (Society
for Electronic Transactions and Security), MGR Knowledge City, CIT
Campus, Taramani, Chennai 600 113, India. }
\email{samy@isichennai.res.in, samy@iitm.ac.in}



\subjclass[2010]{Primary: 30C65,  30C75; Secondary: 30C20, 30C45,
30H10} \keywords{ $K$-quasiconformal harmonic mappings, area, length
distortion.
}

\begin{abstract} In this paper, we investigate the relationships between   linear
measure and harmonic mappings.
\end{abstract}


\maketitle \pagestyle{myheadings} \markboth{ S. Chen, G. Liu  and S.
Ponnusamy }{Linear measure and  $K$-quasiconformal  harmonic
mappings}

\section{Preliminaries and  main results }\label{csw-sec1}

For $a\in\mathbb{C}$ and  $r>0$, we let $\ID(a,r)=\{z:\, |z-a|<r\}$
so that $\mathbb{D}_r:=\mathbb{D}(0,r)$ and thus,
$\mathbb{D}:=\ID_1$ denotes the open unit disk in the complex plane
$\mathbb{C}$. Let $\mathbb{T}=\partial\mathbb{D}$ be the boundary of
$\mathbb{D}$. For a real $2\times2$ matrix $A$, we use the matrix
norm $\|A\|=\sup\{|Az|:\,|z|=1\}$ and the matrix function
$\lambda(A)=\inf\{|Az|:\,|z|=1\}$. For $z=x+iy\in\mathbb{C}$, the
formal derivative of the complex-valued functions $f=u+iv$ is given
by
$$D_{f}=\left(\begin{array}{cccc}
\ds u_{x}\;~~ u_{y}\\[2mm]
\ds v_{x}\;~~ v_{y}
\end{array}\right),
$$
so that
$$\|D_{f}\|=|f_{z}|+|f_{\overline{z}}| ~\mbox{ and }~ \lambda(D_{f})=\big| |f_{z}|-|f_{\overline{z}}|\big |,
$$
where $f_{z}=(1/2)\big( f_x-if_y\big)$ and $f_{\overline{z}}=(1/2)\big(f_x+if_y\big)$ are partial derivatives.

Let $\Omega$ be a domain in $\mathbb{C}$, with non-empty boundary. A sense-preserving
homeomorphism $f$ from a domain $\Omega$ onto $\Omega'$, contained in the Sobolev class $W_{loc}^{1,2}(\Omega)$, is said to be a {\it
$K$-quasiconformal mapping} if, for $z\in\Omega$,
$$\|D_{f}(z)\|^{2}\leq K\det D_{f}(z),~\mbox{i.e.,}~\|D_{f}(z)\|\leq K\lambda\big(D_{f}(z)\big),
$$
where $K\geq1$ and $\det D_{f}$ denotes the determinant of $D_{f}$ (cf. \cite{K,LV,V}).
We note that $\det D_{f} =|f_{z}|^{2}-|f_{\overline{z}}|^{2}$, the Jacobian of $f$, is usually denoted by $J_f$.

A complex-valued function $f$ defined in a  simply connected
subdomain $G$ of $\IC$  is called a {\it harmonic mapping} in $G$ if
and only if both the real and the imaginary parts of $f$ are real
harmonic in $G$. It is well known that every harmonic mapping $f$ in
$G$ admits a decomposition $f=h+\overline{g}$, where $h$ and $g$ are
analytic in $G$. Throughout we use this representation. Without loss
of generality, we assume $0\in G$. If we choose the additive
constant such that $g(0)=0$, then the decomposition is unique.
Because $J_{f} =|h'|^2-|g'|^2$, it follows that $f$ is locally
univalent and sense-preserving in $G$ if and only if
$|g'(z)|<|h'(z)|$ in $G$; or equivalently if $h'(z)\neq0$ and the
dilatation $\omega =g'/h'$ has the property that $|\omega(z)|<1$ in
$G$ (see \cite{Du} and also \cite{Lewy}).


Let $\gamma:\,\varphi(t),~t\in[\alpha,\beta]$, be a curve in
$\mathbb{C}$. Its length $\ell(\gamma)$ is defined by
\be\label{c-eq1}
\ell(\gamma)=\sup\sum_{k=1}^{n}|\varphi(t_{k})-\varphi(t_{k-1})|,
\ee
where the supremum is taken over all partitions
$\alpha=t_{0}<t_{1}<\cdots<t_{n}=\beta$ and all $n\in\{1, 2,\ldots\}$. We call $\gamma$ {\it rectifiable} if
$\ell(\gamma)<+\infty$. It is clear from (\ref{c-eq1}) that
$\diam\gamma\leq\ell(\gamma).$ In the case of a closed curve
$\gamma:\, \psi(\zeta),~\zeta\in\mathbb{T}$, with piecewise
continuously differentiable $\psi$, we can write
$$\ell(\gamma)=\int_{\mathbb{T}}|\psi'(\zeta)|\,|d\zeta|,
$$
where $\mathbb{T}=\partial\mathbb{D}$. If $\gamma_{1},\gamma_{2},\ldots$ are disjoint curves, then we define
$$\ell(\cup_{k}^{\infty}\gamma_{k})=\sum_{k=1}^{\infty}\ell(\gamma_{k}).
$$
In particular, $\ell(\emptyset)=0$ (see \cite[p.~3]{Po1}).

For $p\in(0,\infty]$, the {\it generalized Hardy space
$H^{p}_{g}(\mathbb{D})$} consists of all those functions $f:\,\mathbb{D}\rightarrow\mathbb{C}$ such that $f$ is measurable,
$M_{p}(r,f)$ exists for all $r\in(0,1)$ and  $ \|f\|_{p}<\infty$,
where
$$\|f\|_{p}=
\begin{cases}
\displaystyle\sup_{0<r<1}M_{p}(r,f)
& \mbox{if } p\in(0,\infty)\\
\displaystyle\sup_{z\in\mathbb{D}}|f(z)| &\mbox{if } p=\infty
\end{cases},
~\mbox{ and }~
M_{p}^{p}(r,f)=\frac{1}{2\pi}\int_{0}^{2\pi}|f(re^{i\theta})|^{p}\,d\theta.
$$

\begin{prop}\label{pop-1}
Let $f$ be a $K$-quasiconformal harmonic mapping of
$\mathbb{D}$ onto an inner domain of Jordan curve $\gamma$. Then
$\ell(\gamma)<+\infty~\mbox{if and only if}~\|D_{f}\|\in
H_{g}^{1}(\mathbb{D}).$
\end{prop}
\bpf Assume that $f=h+\overline{g}$, where $h$ and $g$ are analytic
in $\mathbb{D}$. We first prove the necessity. Let
$\ell(\gamma)<+\infty$. Then, by Lemma \Ref{LemA} (in Section \ref{csw-sec2}), we see that
$$F_{n}(z)=\sum_{k=1}^{n}\left|f(ze^{2\pi ki/n})-f(ze^{2\pi (k-1)i/n})\right|
$$
is subharmonic in $\mathbb{D}$ and continuous in $\overline{\mathbb{D}}$. It follows from the maximum principle and
(\ref{c-eq1}) that $F_{n}(z)\leq \ell(\gamma)$ for $z\in\mathbb{D},$ and thus, for $r\in[0,1)$, we have
\begin{eqnarray*}
\frac{r}{K}\int_{0}^{2\pi}\|D_{f}(re^{it})\|\,dt
&\leq&r\int_{0}^{2\pi}\left|h'(re^{it})-e^{-2it}\overline{g'(re^{it})}\right|dt\\
&=&\int_{0}^{2\pi}\left|df(re^{it})\right| =\lim_{n\rightarrow\infty}F_{n}(r)\\
&\leq&\ell(\gamma),
\end{eqnarray*}
which implies that $\|D_{f}\|\in H^{1}_{g}(\mathbb{D}).$

Next, we prove the sufficiency part. For $r\in[0,1)$, let
$\gamma_{r}=\{f(re^{i\theta}): \,\theta\in[0,2\pi)\}$. Since
$\left |zh'(z)-\overline{zg'(z)}\right |$ is subharmonic in $\mathbb{D}$, we see
that
$$\ell(\gamma_{r})=r\int_{0}^{2\pi}\left|h'(re^{it})-e^{-2it}\overline{g'(re^{it})}\right|dt
$$
is an increasing function of $r$ on $[0,1)$. By calculations, we get
$$
\ell(\gamma_{r})\leq\int_{0}^{2\pi}\|D_{f}(re^{it})\|\,dt\leq\sup_{0<r<1}\int_{0}^{2\pi}\|D_{f}(re^{it})\|\,dt<+\infty,
$$
which, together with the monotonicity,  yields that
$\lim_{r\rightarrow1-}\ell(\gamma_{r})$ does exist and thus,
\begin{eqnarray*}
\ell(\gamma)\leq\lim_{r\rightarrow 1-}\ell(\gamma_{r})<+\infty,
\end{eqnarray*}
as desired.
\epf

In \cite{La}, Lavrentiev proved that if $f$ maps $\mathbb{D}$ conformally onto the inner domain of Jordan curve $\gamma$ of finite
length, then, for any $E\subset\gamma$, $\ell(E)>0$ implies that $\ell(f(E))>0$. For univalent harmonic mappings,
we obtain the following result.

\begin{thm}\label{thm-c1}
Let $f$  be a sense-preserving and univalent harmonic mapping from
$\mathbb{\overline{D}}$ onto a domain $\Omega$ and $\partial\Omega$
is a rectifiable Jordan curve. Furthermore, let $E\subset\mathbb{T}$
is measurable with $\ell(E)>0$. Then, we have
$$\ell(f(E))\geq\frac{\ell(\partial\Omega)\ell(E)}{2\pi-\ell(E)}\left[\frac{(|f_{z}(0)|-|f_{\overline{z}}(0)|)(2\pi-\ell(E))}
{\ell(\partial\Omega)}\right]^{\frac{2\pi}{\ell(E)}}.
$$
This estimate is sharp as $\ell(E)\rightarrow2\pi^{-}$ and the
extreme univalent harmonic mapping is
$$f(z)=\frac{M}{2\pi}\int_{0}^{2\pi}\frac{1-|z|^{2}}{|e^{it}-z|^{2}}e^{i\varphi(t)}\,dt
$$
satisfying $|f_{z}(0)|-|f_{\overline{z}}(0)|=M$, where $M$ is a
positive constant and $\varphi(t)$ is a continuously increasing
function in $[0, 2\pi]$ with $ \varphi(2\pi)- \varphi(0)= 2\pi$.
\end{thm}


\begin{cor}
Under the hypotheses of Theorem \ref{thm-c1},  we have $\ell(f(E))>0.$
\end{cor}

Let $\gamma$ be a rectifiable Jordan curve. The shorter arc between
$z$ and $w$ in $\gamma$ will be denoted by $\gamma[z,w]$. We say
that $\gamma$ is a {\it $M$-Lavrentiev curve} if there is a constant
$M>1$ such that $\ell(\gamma[z,w])\leq M|z-w|$ for each $z, w\in\gamma$.
The inner domain of a $M$-Lavrentiev curve is called a {\it
$M$-Lavrentiev domain} (cf. \cite{JK,K,Po1,W}).

The following result is considered to be  a Schwarz-type lemma of subharmonic functions.

\begin{Thm} {\rm  (\cite[Theorem 2]{B})}\label{Thm-cs}
Let $\phi$ be subharmonic in $\mathbb{D}$. If, for all $r\in[0,1)$,
$$A(r)=\sup_{\theta\in[0,2\pi]}\int_{0}^{r}\phi(\rho e^{i\theta})\,d\rho\leq1,
$$
then $A(r)\leq r$.
\end{Thm}

Analogy to Theorem \Ref{Thm-cs}, applying some part of  proof
technique of \cite[Lemma 1]{W}, we obtain a Schwarz type estimate on
the length function of $K$-quasiconformal harmonic  mappings.

\begin{thm}\label{thm-c2}
Suppose that $f$ is a $K$-quasiconformal harmonic mapping of
$\mathbb{D}$ onto a $M$-Lavrentiev domain $\Omega$. Then, for
$r\in(0,2]$, $\rho\in(0,r]$ and any fixed
$\zeta_{0}\in\partial\mathbb{D}$,
$$\int_{0}^{r}\ell(f(\Gamma_{\rho}))\,d\rho
\leq\sqrt{\frac{K\pi A(\Omega)}{3}}\frac{r^{\frac{3}{2}}}{e^{\frac{\alpha}{2}\left(\frac{1}{r}-\frac{1}{2}\right)}}
\leq\sqrt{\frac{K\pi A(\Omega)}{3}}r^{\frac{3}{2}},
$$
where $\alpha=4/[K(1+M)^{2}],$ $A(\Omega)$ is the area of $\Omega$
and $\Gamma_{\rho}$ is the arc of the circle
$\partial\mathbb{D}(\zeta_{0},\rho)$ which lies in
$\overline{\mathbb{D}}$.
\end{thm}

We say that a simply connected domain $G\subset\mathbb{C}$ is {\it
$M$-linearly connected} if, for any two points $z_{1}, z_{2}\in G$,
there is a curve $\gamma\subset G$ and a constant $M\geq1$ such that
\be\label{eqc-24}
\diam \gamma\leq M|z_{1}- z_{2}|.
\ee


Zinsmeister \cite{Z} obtained an analytic characterization of $M$-Lavrentiev domains (see also \cite[Chapter 7]{Po1}). The
following result is an analogous result to the analytic characterization of $M$-Lavrentiev domains.


\begin{thm}\label{thm-c3}
Let $f$ be a $K$-quasiconformal harmonic mapping from $\mathbb{D}$
onto the inner domain $G$ of a rectifiable Jordan curve. If $G$ is a
$M_{1}$-Lavrentiev domain and for each $\zeta\in\partial\mathbb{D},$
\be\label{eq-h} \|D_{f}(\rho\zeta)\|\leq M_{1}'
\|D_{f}(r\zeta)\|\left(\frac{1-\rho}{1-r}\right)^{\delta-1}~(0\leq
r\leq\rho<1), \ee then $G$ is $M_{2}$-linearly connected and, for
all $z\in\mathbb{D},$ there is a constant $M'$ such that
\be\label{eqc-26}
\frac{1}{\ell(I(z))}\int_{I(z)}\|D_{f}(\zeta)\|\,|d\zeta|\leq
M'\|D_{f}(z)\|, \ee where $I(z)=\{\zeta\in\mathbb{T}: \,
|\arg\zeta-\arg z|\leq\pi(1-|z|)\}$ and $\delta\in(0,1)$, $M_{1}'$,
$M_{1}$, $M_{2}$ are constants.
\end{thm}

For any fixed $\theta\in[0,2\pi]$,  the {\it radial length} of the
curve $C_{\theta}(r)=\big\{w=f(\rho e^{i\theta}):\, 0\leq\rho\leq
r\big\}$ with counting multiplicity is defined by
$$\ell_{f}^{\ast}(\theta,r)=\int_{0}^{r}\left |\,df(\rho
e^{i\theta})\right |=\int_{0}^{r}\left |f_{z}(\rho
e^{i\theta})+e^{-2i\theta}f_{\overline{z}}(\rho
e^{i\theta})\right |\,d\rho,
$$
where $r\in[0,1)$ and $f$ is a harmonic mapping defined in
$\mathbb{D}$. In particular, let
$$\ell_{f}^{\ast}(\theta,1)=\sup_{0<r<1}\ell_{f}^{\ast}(\theta,r).
$$

\begin{prop}\label{pop-2}
Suppose that $f$ is a bounded harmonic mapping in $\mathbb{D}$ and
$r_{0}\in(0,1)$. For $\zeta\in\mathbb{D},$ let
$F(\zeta)=f(r_{0}\zeta)$.  Then, for $\rho\in[0,1)$ and $\theta\in[0,2\pi]$,
$$\ell_{F}^{\ast}(\theta,r)=\int_{0}^{r}r_{0}\left |f_{z}(\rho r_{0}
e^{i\theta})+e^{-2i\theta}f_{\overline{z}}(\rho r_{0}
e^{i\theta})\right |\,d\rho\leq Mr,
$$
where $r\in(0,1)$ and
$M=\frac{2r_{0}}{\pi}\sup_{z\in\mathbb{D}}|f(z)|\log((1+r_{0})/(1-r_{0})).$
\end{prop}
\bpf
For $\zeta\in\mathbb{D}$, let $F(\zeta)=f(r_{0}\zeta)$. 
By \cite[Theorem 3]{Co}, we get
\begin{eqnarray*}
\ell_{F}^{\ast}(\theta,r)&=&\int_{0}^{r}r_{0}\big |f_{z}(\rho r_{0}
e^{i\theta})+e^{-2i\theta}f_{\overline{z}}(\rho r_{0}
e^{i\theta})\big |\,d\rho\\
&\leq&\int_{0}^{r}r_{0}\|D_{f}(r_{0}\rho e^{i\theta})\|\,d\rho\\
&\leq&\frac{4r_{0}\sup_{z\in\mathbb{D}}|f(z)|}{\pi}\int_{0}^{r}\frac{d\rho}{1-r_{0}^{2}\rho^{2}}\\
&=&\frac{2r_{0}\sup_{z\in\mathbb{D}}|f(z)|}{\pi}\log\left(\frac{1+r_{0}r}{1-r_{0}r}\right)\\
&\leq&\frac{2r_{0}\sup_{z\in\mathbb{D}}|f(z)|}{\pi}\log\left(\frac{1+r_{0}}{1-r_{0}}\right)
=M.
\end{eqnarray*} 
By the subharmonicity of $\|D_{f}(r_{0}\rho e^{i\theta})\|$ and
Theorem \Ref{Thm-cs}, we see that, for $\theta\in[0,2\pi]$,
$\ell_{F}^{\ast}(\theta,r)\leq Mr,$ where $r\in(0,1)$.
\epf

In \cite{CPR}, the authors obtained the coefficient estimates of a
class of $K$-quasiconformal harmonic mapping with a finite perimeter
length. In the following, we will investigate the coefficient
estimates on a class of $K$-quasiconformal harmonic mapping with the
finite radial length.

\begin{thm}\label{thm-5}
Let
$f(z)=\sum_{n=0}^{\infty}a_{n}z^{n}+\sum_{n=1}\overline{b_{n}}\overline{z}^{n}$
be a $K$-quasiconformal harmonic mapping on $\ID$. If, for all
$\theta\in[0,2\pi]$, $\ell_{f}^{\ast}(\theta,1)\leq M$ for some
positive constant $M$, then \be\label{eq-cs} |a_{n}|+|b_{n}|\leq KM
~\mbox{ for $n\geq 1$,} \ee In particular, if $K=1$, then the
estimate \eqref{eq-cs} is sharp and the extreme function is
$f(z)=Mz$.
\end{thm}

Let $d_{\Omega}(z)$ be the Euclidean distance from $z$ to the
boundary $\partial \Omega$ of the domain $\Omega$. In the following,
we investigate the behaviour   on the ratio of the radial length and
the perimeter length on
 $K$-quasiconformal harmonic mappings.

\begin{thm}\label{thm-c4}
Let $f$ be a $K$-quasiconformal harmonic mapping of $\mathbb{D}$ onto a bounded domain $D$. Then, for
$r\in(0,1)$ and all $\theta\in[0,2\pi]$,
$$\frac{\ell_{f}^{\ast}(\theta,r)}{\ell(r)}\leq\frac{32r(1+r)K^{3}\sup_{z\in\mathbb{D}}
|f(z)|}{\int_{0}^{2\pi}d_{D}(f(r e^{it}))\,dt}
$$
and
$$\lim_{r\rightarrow0^{+}}\left\{\frac{\sup_{\theta\in[0,2\pi]}\ell_{f}^{\ast}(\theta,r)}{\ell(r)}\right\}=0.
$$
\end{thm}

The proofs of Theorems \ref{thm-c1}, \ref{thm-c2}, \ref{thm-c3},
\ref{thm-5} and \ref{thm-c4}  will be presented in Section
\ref{csw-sec2}.

\section{The proofs of the main results }\label{csw-sec2}
The following lemmas are well-known.

\begin{Lem}{\rm (see \cite{KMM}, \cite[Proposition 1.8]{K} and \cite{Pav})}\label{LemA}
If $f=h+\overline{g}$ is a $K$-quasiconformal harmonic mapping of
$\mathbb{D}$ onto a Jordan domain with a rectifiable boundary, then
$h$ and $g$ have absolutely continuous extension to $\mathbb{T}$.
\end{Lem}

\begin{Lem}{\rm (cf. \cite{Ru})}\label{Jensen}
Let $(\Omega,A,\mu)$ be a measure space such that $\mu(\Omega)=1.$
If $g$ is a real-valued function that is $\mu$-integrable, and if
$\varphi$ is a convex function on the real line, then
$$\varphi\left (\int_{\Omega}g\,d\mu\right )\leq\int_{\Omega}\varphi\circ g\,d\mu.
$$
\end{Lem}

\begin{Lem}{\rm (cf. \cite{Ca})}\label{Lem-B}
Among all rectifiable Jordan curves of a given length, the circle
has the maximum interior area.
\end{Lem}

\subsection*{Proof of Theorem \ref{thm-c1}}
Let $f=h+\overline{g}$  be a sense-preserving and univalent harmonic
mapping in $\mathbb{\overline{D}}$, where $h$ and $g$ are analytic
in $\mathbb{\overline{D}}$. Then, $h'(z)\neq 0$ for
$z\in\mathbb{D}$ 
and the dilatation $\omega$ defined by $\omega=g'/h'$ is analytic and $|\omega (z)|<1$ in $\ID$. Since
$\log|h'|$ is harmonic in $\mathbb{D}$ and $\log(1-|\omega|)$ is
subharmonic in $\mathbb{D}$, we see that, for $r\in(0,1)$,
$$\log|h'(0)|=\frac{1}{2\pi}\int_{E}\log|h'(rz)|\,|dz|+\frac{1}{2\pi}\int_{\mathbb{T}\backslash E}\log|h'(rz)|\,|dz|
$$
and
$$\log(1-|\omega(0)|)\leq\frac{1}{2\pi}\int_{E}\log(1-|\omega(rz)|)\,|dz|+\frac{1}{2\pi}\int_{\mathbb{T}\backslash E}\log(1-|\omega(rz)|)\,|dz|,
$$
which, together with Lemma \Ref{Jensen} (Jensen's inequality), imply
that

\vspace{8pt}

$2\pi\log\big[|h'(0)|(1-|\omega(0)|)\big]$
\begin{eqnarray*}
&\leq&\int_{E}\log\big[|h'(rz)|(1-|\omega(rz)|)\big]\,|dz|
+\int_{\mathbb{T}\backslash E}\log\big[|h'(rz)|(1-|\omega(rz)|)\big]\,|dz|\\
&\leq&\ell(E)\log\left[\frac{1}{\ell(E)}\int_{E}\big(|h'(rz)|-|g'(rz)|\big)\,|dz|\right]\\
&&+\big(2\pi-\ell(E)\big)\log\left[\frac{1}{2\pi-\ell(E)}\int_{\mathbb{T}\backslash
E}\big(|h'(rz)|-|g'(rz)|\big)\,|dz|\right]\\
&\leq&\ell(E)\log\left[\frac{1}{r\ell(E)}\int_{E}r\big|h'(rz)-\overline{g'(rz)}\overline{z}^{2}\big|\,|dz|\right]\\
&&+\big(2\pi-\ell(E)\big)\log\left[\frac{1}{r[2\pi-\ell(E)]}\int_{\mathbb{T}\backslash
E}r\big|h'(rz)-\overline{g'(rz)}\overline{z}^{2}\big|\,|dz|\right].
\end{eqnarray*}
By letting $r\rightarrow1^{-}$ and Lemma \Ref{LemA}, we have
\begin{eqnarray*}
2\pi\log\big[|h'(0)|(1-|\omega(0)|)\big]&\leq&\ell(E)\log\frac{\ell(f(E))}{\ell(E)}+(2\pi-\ell(E))\log\frac{\ell(f(\mathbb{T}\backslash
E))}{(2\pi-\ell(E))}\\
&\leq&
\ell(E)\log\ell(f(E))-\ell(E)\log\ell(E)\\
&& -(2\pi-\ell(E))\log(2\pi-\ell(E))+(2\pi-\ell(E))\log L,
\end{eqnarray*}
which implies  that
\be\label{eq-CLP-1}
\ell(f(E))\geq\frac{L\ell(E)}{2\pi-\ell(E)}\left[\frac{(|h'(0)|-|g'(0)|)(2\pi-\ell(E))}{L}\right]^{\frac{2\pi}{\ell(E)}},
\ee
where $\ell(\partial\Omega)=L$.
Now we prove the sharpness part.  By calculation, we have
$$\lim_{\ell(E)\rightarrow2\pi^{-}}\frac{L\ell(E)}{2\pi-\ell(E)}\left[\frac{(|h'(0)|-|g'(0)|)(2\pi-\ell(E))}{L}\right]^{\frac{2\pi}{\ell(E)}}=2\pi(|h'(0)|-|g'(0)|).
$$
Let \be\label{eq-CLP-2}
f(z)=\frac{M}{2\pi}\int_{0}^{2\pi}\frac{1-|z|^{2}}{|e^{it}-z|^{2}}e^{i\varphi(t)}\,dt
\ee satisfying $|f_{z}(0)|-|f_{\overline{z}}(0)|=M$, where $M$ is a
positive constant and $\varphi(t)$ is a  continuously increasing
function in $[0, 2\pi]$ with $ \varphi(2\pi)- \varphi(0)= 2\pi$. If
$\ell(E)\rightarrow2\pi^{-}$, then the  sense-preserving and
univalent harmonic mapping  $f$ of (\ref{eq-CLP-2}) shows that the
estimate of (\ref{eq-CLP-1}) is sharp. \hfill $\Box$

\subsection*{Proof of Theorem \ref{thm-c2}} Let $f=h+\overline{g}$ satisfy the assumption, where $h$ and $g$ are analytic in $\mathbb{D}$.
 By Lemma \Ref{LemA}, we know that $h$ and $g$  can be absolutely
continuous extension to $\mathbb{T}$. For $r\in(0,2]$ and
$\rho\in(0,r],$ let $\Delta_{\rho}=\{z:
\,|z-\zeta_{0}|\leq\rho~\mbox{and}~|z|\leq1\}.$ Let $\Gamma_{\rho}$
denote the arc of the circle $\partial\mathbb{D}(\zeta_{0},\rho)$
which lies in $\overline{\mathbb{D}}$. Then we have
\begin{eqnarray*}
\ell^{2}\big(f(\Gamma_{\rho})\big)&=&\left(\int_{\Gamma_{\rho}}
\big|f_{z}(\zeta_{0}+\rho e^{it})-e^{-2it}f_{\overline{z}}(\zeta_{0}+\rho
e^{it})\big|\rho \,dt\right)^{2}\\
&\leq&\left(\int_{\Gamma_{\rho}}\rho
dt\right)\left(\int_{\Gamma_{\rho}} \big|f_{z}(\zeta_{0}+\rho
e^{it})-e^{-2it}f_{\overline{z}}(\zeta_{0}+\rho e^{it})\big|^{2}\rho
\,dt\right)\\
&=&\ell(\Gamma_{\rho}) \rho\int_{\Gamma_{\rho}}
\big|f_{z}(\zeta_{0}+\rho e^{it})-e^{-2it}f_{\overline{z}}(\zeta_{0}+\rho e^{it})\big|^{2}\rho\,dt\\
&\leq&2\rho^{2}\arccos\frac{\rho}{2} \int_{\Gamma_{\rho}} \|D_{f}(\zeta_{0}+\rho e^{it})\|^{2}\rho \,dt\\
&\leq&K\pi\rho^{2}\int_{\Gamma_{\rho}} J_{f}(\zeta_{0}+\rho e^{it})\rho \,dt,
\end{eqnarray*}
which implies that
\beq\label{eqc-4}
P(r)=\int_{0}^{r}\frac{\ell^{2}\big(f(\Gamma_{\rho})\big)}{\rho^{2}}d\rho&\leq&K\pi\int_{0}^{r}\int_{\Gamma_{\rho}}
J_{f}(\zeta_{0}+\rho e^{it})\rho \,dt\,d\rho\\ \nonumber
&\leq&K\pi A_{f}(r),
\eeq
where $z=\zeta_{0}+\rho e^{it}$
and  $A_{f}(r)$ denotes the area of $f(\Delta_{r})$. Since the boundary of $\Omega$ is a $M$-Lavrentiev curve, we see that
$$\ell\left(f(\partial\Delta_{r})\right)\leq\ell\big(f(\Gamma_{r})\big)+M\ell\big(f(\Gamma_{r})\big) =(1+M)\ell\big(f(\Gamma_{r})\big)
$$
and thus, by Lemma \Ref{Lem-B}, we have
\beq\label{eqc-5}
A_{f}(r)\leq  
\frac{(1+M)^{2}}{4\pi}\ell^{2}\big(f(\Gamma_{r})\big).
\eeq
By (\ref{eqc-4}) and (\ref{eqc-5}), we obtain
$$P(r)\leq\frac{K(1+M)^{2}}{4}\ell^{2}\big(f(\Gamma_{r})\big).
$$
By calculations, for $r\in(0,2]$, we have
$$\frac{4}{K(1+M)^{2}}P(r)\leq\ell^{2}\big(f(\Gamma_{r})\big)=r^{2}P'(r),
$$
which gives that
\beq\label{eqc-6}
\frac{\alpha}{r^{2}}\leq\frac{P'(r)}{P(r)},
\eeq
where $\alpha=4/[K(1+M)^{2}].$ By (\ref{eqc-6}), we get
$$\int_{r}^{2}\frac{\alpha}{\rho^{2}}d\rho\leq\int_{r}^{2}\frac{P'(\rho)}{P(\rho)}\,d\rho,
$$
which by integration, together with (\ref{eqc-4}), yield that
\beq\label{eqc-7}
P(r)\leq\frac{P(2)}{e^{\alpha\left(\frac{1}{r}-\frac{1}{2}\right)}}
\leq\frac{K\pi A(\Omega)}{e^{\alpha\left(\frac{1}{r}-\frac{1}{2}\right)}},
\eeq
where $A(\Omega)$ is the area of $\Omega$. 

By Cauchy-Schwarz's inequality and (\ref{eqc-7}), we have
$$\left(\int_{0}^{r}\ell(f(\Gamma_{\rho}))d\rho\right)^{2}\leq P(r)\int_{0}^{r}\rho^{2}\, d\rho
\leq\frac{K\pi A(\Omega)r^{3}}{3e^{\alpha\left(\frac{1}{r}-\frac{1}{2}\right)}},
$$
which gives that
$$\int_{0}^{r}\ell(f(\Gamma_{\rho}))\,d\rho\leq\sqrt{\frac{K\pi
A(\Omega)}{3}}\frac{r^{\frac{3}{2}}}{e^{\frac{\alpha}{2}\left(\frac{1}{r}-\frac{1}{2}\right)}}\leq\sqrt{\frac{K\pi
A(\Omega)}{3}}r^{\frac{3}{2}}.
$$
The proof of this theorem is complete. \hfill $\Box$

\vspace{8pt}

A Jordan curve $J$ is said to be a {\it $M$-quasicircle} if, for any $z_{1}, z_{2}\in J$, there is a constant $M\geq1$ such that
\be\label{eqc-25}
\diam J(z_{1}, z_{2})\leq M|z_{1}- z_{2}|.
\ee
The inner domain $G$ of a quasicircle $J$ is called a {\it
$M$-quasidisk}. We say that the curve $\gamma\subset\mathbb{C}$ is
{\it Ahlfors-regular} if there is a positive constant $M$ such that
$$\ell(\gamma\cap\mathbb{D}(w,r))\leq Mr,
$$
where $r\in(0,\infty)$ and $w\in\mathbb{C}$ (cf. \cite{Po1}).

\begin{Lem}{\rm  (\cite[Proposition 7.7]{Po1})}\label{Lem-C}
 A domain is a $M_{3}$-Lavrentiev domain if and only if it is an
 $M_{4}$-Ahlfors-regular quasidisk, where $M_{3}$ and $M_{4}$ are positive constants.
\end{Lem}

\subsection*{Proof of Theorem \ref{thm-c3}}
Assume that $f=h+\overline{g}$ is a $K$-quasiconformal harmonic mapping from
$\mathbb{D}$ onto the inner domain $G$ of a rectifiable Jordan
curve, where $h$ and $g$ are analytic in $\mathbb{D}$.
Let $w_{1}, w_{2}\in G$ and let $[w_{1}, w_{2}]$ denote the
line segment with endpoints $w_{1}$ and $w_{2}$. If $[w_{1}, w_{2}]\subset G$, then (\ref{eqc-24}) holds. Without loss of
generality, we assume that $[w_{1}, w_{2}]\nsubseteq G.$ Let $f(\xi_{k})$ be the boundary point on $[w_{1}, w_{2}]$ nearest to
$w_{k}$, where $k=1,2$ and $\xi_{k}\in\mathbb{T}.$ By (\ref{eqc-25}), we see that one of the arcs
$\mathbb{T}_{[\xi_{1},\xi_{2}]}$ of $\mathbb{T}$ from $\xi_{1}$ to
$\xi_{2}$ satisfies
$$\diam f(\mathbb{T}_{[\xi_{1},\xi_{2}]})\leq M_{1}|f(\xi_{1})-f(\xi_{2})|.
$$
If $r$ is close enough to $1$, then $f(r\mathbb{T}_{[\xi_{1},\xi_{2}]})$ is a curve in $G$ of
$$\diam f(r\mathbb{T}_{[\xi_{1},\xi_{2}]})\leq2M_{1}|f(\xi_{1})-f(\xi_{2})|<2M_{1}|w_{1}- w_{2}|
$$
which can be connected within $G$ by curves $\gamma_{k}$ to $w_{k}$ satisfying
$$\ell(\gamma_{k})<|w_{1}- w_{2}|.
$$
Then $\gamma=\gamma_{1}\cup f(r\mathbb{T}_{[\xi_{1},\xi_{2}]})\cup\gamma_{2}$ is curve in $G$
from $w_{1}$ to $w_{2}$ satisfying
$$\diam\gamma\leq(2M_{1}+2)|w_{1}- w_{2}|,
$$
which implies that $G$ is a $(2M_{1}+2)$-linearly connected domain.

Now we  prove (\ref{eqc-26}). Let  $I(z)=\{\zeta\in\mathbb{T}:
\,|\arg\zeta-\arg z|\leq\pi(1-|z|)\}.$ By (\ref{eq-h}) and
\cite[Theorems 1 and 2]{CP}, we see that there is a constant
$M^{\ast}$ such that \be\label{eqc-27} \diam f(I(z))\leq
M^{\ast}d_{G}(z). \ee Applying the inequality (2.3) in \cite{CP}, we
get \be\label{eqc-28}
d_{G}(z)\leq\frac{2K}{1+K}\|D_{f}(z)\|(1-|z|^{2}). \ee It follows
from (\ref{eqc-27}) and (\ref{eqc-28}) that \be\label{eqc-29} \diam
f(I(z))\leq \frac{2KM^{\ast}}{1+K}\|D_{f}(z)\|(1-|z|^{2})=r_{z}, \ee
which implies that $f(I(z))$ lies in $\mathbb{D}_{r_{z}}$. By
(\ref{eqc-29}) and  Lemma \Ref{Lem-C}, we conclude that there is a
constant $M'$ such that
\begin{eqnarray*}
\frac{1}{K}\int_{I(z)}\|D_{f}(\zeta)\|\,|d\zeta|
&\leq&\int_{I(z)}\left|h'(\zeta)- \overline{\zeta}^{2}\overline{g'(\zeta)}\right|\,|d\zeta|=\ell(f(I(z)))\\
&\leq& M'\|D_{f}(z)\|(1-|z|^{2})\\
&\leq& M'\ell(I(z))\|D_{f}(z)\|.
\end{eqnarray*}
The proof of the  theorem is complete.  \hfill $\Box$


\begin{Thm} {\rm (\cite[Proposition 3.1]{MM} and \cite[Theorem 3.2]{MM})  }\label{Thm-MM}
Let $f$ be a $K$-quasiconformal harmonic mapping from $\mathbb{D}$
onto itself. Then for all $z\in\mathbb{D}$, we have
$$\frac{1+K}{2K}\left (\frac{1-|f(z)|^{2}}{1-|z|^{2}}\right )\leq|f_{z}(z)|\leq\frac{K+1}{2}\left (\frac{1-|f(z)|^{2}}{1-|z|^{2}}\right ).
$$
\end{Thm}

\subsection*{Proof of Theorem \ref{thm-5}}
Let $f(z)=\sum_{n=0}^{\infty}a_{n}z^{n}+\sum_{n=1}\overline{b}_{n}\overline{z}^{n}$
be a $K$-quasiconformal harmonic mapping on $\ID$. Then, by Cauchy's integral formula, for $\rho \in (0,1)$ and $n\geq 1$, we get
$$na_{n}=\frac{1}{2\pi i}\int_{|z|=\rho}\frac{f_{z}(z)}{z^{n}}\,dz~\mbox{ and }~
nb_{n}=\frac{1}{2\pi i}\int_{|z|=\rho}\frac{\overline{f_{\overline{z}}(z)}}{z^{n}}\,dz,
$$
which imply that
\beq\label{eqc-30}
n(|a_{n}|+|b_{n}|)&=&\frac{1}{2\pi}\left|\int_{|z|=\rho}\frac{f_{z}(z)}{z^{n}}\,dz\right|
+\frac{1}{2\pi}\left|\int_{|z|=\rho}\frac{\overline{f_{\overline{z}}(z)}}{z^{n}}\,dz\right|\\ \nonumber
&\leq&\frac{1}{2\pi \rho^{n-1}}\int_{0}^{2\pi}\|D_{f}(\rho e^{i\theta})\|\,d\theta.
\eeq
By calculations, for $\theta\in[0,2\pi]$, we obtain
\begin{eqnarray*}
\ell_{f}^{\ast}(\theta,r)
&=&\int_{0}^{r}\left |f_{z}(\rho e^{i\theta})+e^{-2i\theta}f_{\overline{z}}(\rho e^{i\theta})\right |\,d\rho\\
&\geq&\int_{0}^{r}\lambda(D_f)(\rho e^{i\theta})\,d\rho\\
&\geq&\frac{1}{K}\int_{0}^{r}\|D_{f}(\rho e^{i\theta})\|\,d\rho,
\end{eqnarray*}
which gives
\be\label{eqc-31}
\int_{0}^{r}\|D_{f}(\rho e^{i\theta})\|\,d\rho\leq K\ell_{f}^{\ast}(\theta,r)\leq KM.
\ee
By (\ref{eqc-31}), the subharmonicity of $D_{f}(\rho e^{i\theta})$ and Theorem \Ref{Thm-cs}, we have
\be\label{eqc-32}
\int_{0}^{r}\|D_{f}(\rho e^{i\theta})\| \, d\rho\leq KMr.
\ee
By (\ref{eqc-30}) and (\ref{eqc-32}), we get
\begin{eqnarray*}
2\pi n(|a_{n}|+|b_{n}|)\int_{0}^{r}\rho^{n-1}\, d\rho
&=&\int_{0}^{r}\left(\int_{0}^{2\pi}\|D_{f}(\rho e^{i\theta})\|d\theta\right)\,d\rho\\
&=&\int_{0}^{2\pi}\left(\int_{0}^{r}\|D_{f}(\rho e^{i\theta})\|d\rho\right) \, d\theta\\
&\leq&2\pi KMr,
\end{eqnarray*}
which yields that
$$|a_{n}|+|b_{n}|\leq\frac{KM}{r^{n-1}} ~\mbox{ for $n\geq  1$}.
$$
Since this is true for each $r<1$, the desired bound follows by
letting $r\rightarrow 1^{-}$. \hfill $\Box$

\subsection*{Proof of Theorem \ref{thm-c4}}
Let $f=h+\overline{g}$ be a $K$-quasiconformal harmonic mapping of
$\mathbb{D}$ onto a bounded domain $D$,  where $h$ and $g$ are analytic in $\mathbb{D}.$
 By \cite[Proposition 13]{Mi} and Theorem \Ref{Thm-MM}, for
$r\in(0,1)$ and all $\theta\in[0,2\pi]$, we obtain
\beq\label{eqc-20}
\ell_{f}^{\ast}(\theta,r)
&=&\int_{0}^{r}|df(\rho e^{i\theta})|=\int_{0}^{r}\left |h'(\rho
e^{i\theta})+e^{-2i\theta}\overline{g'(\rho e^{i\theta})}\right |\,d\rho\\
\nonumber &\leq&\int_{0}^{r}\|D_{f}(\rho e^{i\theta})\|\,d\rho\\
\nonumber &\leq&16K\int_{0}^{r}\frac{d_{D}(f(\rho e^{i\theta}))}{1-\rho^{2}}\,d\rho\\
\nonumber &\leq&16K\sup_{z\in\mathbb{D}}|f(z)|\int_{0}^{r}\frac{1}{1-\rho^{2}}\,d\rho\\
\nonumber &=&8K\sup_{z\in\mathbb{D}}|f(z)|\log\frac{1+r}{1-r}
\eeq
and
\beq\label{eqc-21}
\ell(r)&=&\int_{0}^{2\pi}r\left |h'(r e^{i\theta})-e^{-2it}\overline{g'(r e^{it})}\right |\,dt\\
\nonumber &\geq&r\int_{0}^{2\pi}\lambda(D_f)(re^{it})\,dt\\
\nonumber &\geq&\frac{r}{K}\int_{0}^{2\pi}\|D_{f}(re^{it})\|\,dt\\
\nonumber &\geq& \frac{r(1+K)}{2K^{2}(1-r^{2})}\int_{0}^{2\pi}d_{D}(f(re^{it}))\,dt
\eeq
where the last inequality is a consequence of \cite[Inequality (2.3)]{CP}.
Equations (\ref{eqc-20}) and (\ref{eqc-21}) imply that
$$
\frac{\ell_{f}^{\ast}(\theta,r)}{\ell(r)}\leq16K^{3}\sup_{z\in\mathbb{D}}
|f(z)|\frac{(1-r^{2})\log\frac{1+r}{1-r}}{\int_{0}^{2\pi}d_{D}(f(r
e^{it}))\,dt}\leq\frac{32r(1+r)K^{3}\sup_{z\in\mathbb{D}}
|f(z)|}{\int_{0}^{2\pi}d_{D}(f(r e^{it}))\,dt}
$$
and, for all $\theta\in[0,2\pi]$,
$$\lim_{r\rightarrow0^{+}}\frac{\ell_{f}^{\ast}(\theta,r)}{\ell(r)}=0.
$$
The proof of this theorem is complete. \hfill $\Box$




\bigskip

{\bf Acknowledgements:}  This research was partly supported by the
National Natural Science Foundation of China ( No. 11571216 and No.
11401184), the Hunan Province Natural Science Foundation of China
(No. 2015JJ3025), the Excellent Doctoral Dissertation of Special
Foundation of Hunan Province (higher education 2050205), the
Construct Program of the Key Discipline in Hunan Province. The third
is on leave from IIT Madras. The second author thanks Indian Statistical Institute, Chennai Centre for the
hospitality and the support during the period of my visit to India.



\normalsize


\begin{thebibliography}{99}

\bibitem{B} {\sc E. F. Beckenbach,}
A relative of the lemma of Schwarz,
\textit{Bull. Amer. Math. Soc.,} {\bf 44} (1938), 698--707.

\bibitem{Ca} {\sc T. Carleman,}
Zur Theorie der Minimalfl\"achen,
\textit{Math. Z.,} {\bf 9} (1921), 154--160.

\bibitem{CPR} {\sc S. Chen, S. Ponnusamy and A. Rasila,}
Lengths, areas and Lipschitz-type spaces of planar harmonic mappings,
\textit{Nonlinear Anal.,} {\bf 115} (2015), 62--70.


\bibitem{CP} {\sc S. Chen and S. Ponnusamy,}
John disks and $K$-quasiconformal harmonic mappings, \textit{J.
Geom. Anal.,} 2016, DOI: 10.1007/s12220-016-9727-6, published
online.



\bibitem{Co}{\sc F. Colonna,}
The Bloch constant of bounded harmonic mappings,
\textit{Indiana Univ. Math. J.}, {\bf 38} (1989), 829--840.

\bibitem{Du} {\sc P. Duren,}
{\it Harmonic mappings in the plane,} Cambridge Univ. Press, 2004.

\bibitem{JK} {\sc D. S. Jerison and C. E. Kenig,}
Hardy spaces, $A_{\infty}$ and singular integrals in chord-arc domains,
\textit{Math. Scand.,} {\bf 50} (1982), 221--247.

\bibitem{KMM} {\sc D. Kalaj, M. Markovi\'c and M. Mateljevi\'c,}
Carath\'eodory and Simirov type theorems for harmonic mappings of the unit disk onto surfaces,
\textit{Ann. Acad. Sci. Fenn. Math.,} {\bf 38} (2013), 565--580.

\bibitem{K} {\sc D. Kalaj,}
Muckenhoupt weights and Lindel\"of theorem for harmonic mappings,
\textit{Adv. Math.,} {\bf 280} (2015), 301--321.

\bibitem{MM} {\sc M. Kne{\rm $\breve{z}$}evi\'c and M. Mateljevi\'c,}
On the quasi-isometries of harmonic quasiconformal mappings,
\textit{J. Math. Anal. Appl.,} {\bf 334} (2007), 404--413.

\bibitem{La} {\sc M. Lavrentiev,}
Boundary problems in the theory of univalent functions,
\textit{Amer. Math. Soc. Tansl. Ser.,} {\bf 32} (1963), 1--35.

\bibitem{LV} {\sc O. Lehto and K. I. Virtanen,}
{\it Quasiconformal mappings in the plane}, Springer Verlag, 1973.

\bibitem{Lewy} {\sc H. Lewy,}
On the non-vanishing of the Jacobian in certain one-to-one mappings,
\textit{Bull. Amer. Math. Soc.,} {\bf 42} (1936), 689--692.

\bibitem{Mi} {\sc M. Mateljevi\'c,}
Distrotion of quasiregular mappings and equivalent norms on Lipschitz-type spaces,
\textit{Abstr. Appl. Anal.,} Volume 2014, Article ID 895074, 20 pages.


\bibitem{Pav} {\sc M. Pavlovi\'c,}
Boundary correspondence under harmonic quasiconformal homeomorphisms of the unit disc,
\textit{Ann. Acad. Sci. Fenn. Math.,} {\bf 27} (2002), 365--372.


\bibitem{Po1} {\sc Ch. Pommerenke,}
{\it Boundary behaviour of conformal maps}, Springer-Verlag, 1992.

\bibitem{Ru} {\sc W. Rudin,}  Real and Complex Analysis. McGraw-Hill. ISBN
0-07-054234-1,1987.


\bibitem{V} {\sc M. Vuorinen,}
{\it Conformal geometry and quasiregular mappings,}
\textit{Lecture Notes in Math.} Vol. 1319, Springer-Verlag, 1988.

\bibitem{W} {\sc S. E. Warschawski,}
On differentiability at the boundary in conformal mapping,
\textit{Proc. Amer. Math. Soc.,} {\bf 12} (1961), 614--620.


\bibitem{Z} {\sc M. Zinsmeister,}
\textit{ Domaines de Lavrentiev},
Publications Math. Orsay, 162pp, no.3, 1985.

\end{thebibliography}
\end{document}